\documentclass[11pt, leqno, twoside]{amsart}
\usepackage{amssymb, amsmath, amsthm} 
\usepackage{latexsym} 
\usepackage{amscd}

 \newtheorem{thm}{Theorem}[section]
 
 \newtheorem{lem}[thm]{Lemma}
 
 \theoremstyle{definition}
 
 \theoremstyle{remark}
 \newtheorem{rem}[thm]{Remark}
 
 \numberwithin{equation}{section}

\begin{document}

\title[Primes represented by quadratic polynomials]
 {Primes represented by quadratic\\ polynomials via exceptional characters}

\author{Fernando Chamizo}

\address{Departamento de Matem\'aticas and ICMAT\\
Universidad Aut\'onoma de Ma\-drid\\
28049 Madrid, Spain}

\email{fernando.chamizo@uam.es}

\thanks{The first author is partially supported by the MTM2017-83496-P grant of the MICINN (Spain) and
by ``Severo Ochoa Programme for Centres of Excellence in R\&{D}'' (SEV-2015-0554). This latter grant supported the visit of the second author to the ICMAT where this work was completed. The second author is partially supported by the PID2019-110224RB-I00 grant of the MICINN (Spain)}

\author{Jorge Jim\'enez Urroz}
\address{Departamento de Matem\'aticas \\
Universitat Polit\`ecnica Catalunya\\
Barce\-lona, Spain}
\email{jorge.urroz@upc.edu}

\subjclass[2010]{Primary 11N32; Secondary 11N35, 11M20}

\keywords{Exceptional characters, $L$ functions, Sieve methods}

\begin{abstract}
We estimate the number of primes represented by a general quadratic polynomial with discriminant  $\Delta$, assuming that the corresponding real character is exceptional.
\end{abstract}

\maketitle

\section{Introduction}

Let $f(x)=ax^2+bx+c$ be a quadratic polynomial with integer coefficients such that $(a,b,c)=1$, $a+b$ or $c$ odd and discriminant  $\Delta\not=\square$. Conjecture~F in the classic  work \cite{HaLi} claims that there are infinitely many prime numbers of the form $f(n)$ when $a>0$. 
Note that it is elementary that 
the imposed conditions on $f$ are necessary to represent infinitely many primes. 
If $a<0$, under the same conditions, we still expect to capture many primes if $\{n\in\mathbb{Z}\,:\, f(n)>0\}$ is large and  Conjecture~H in \cite{HaLi} is an instance of it.
\smallskip

In \cite{GrMo} and in \cite{FrIw} special cases of these conjectures are addressed assuming the existence of exceptional characters. For instance, in the second paper it is proved that positive exceptional fundamental discriminants $D$ can be written as $D=m^2+p$ and that if $D$ is ``exceptional enough'' we have an asymptotic formula for the number of primes of the form $D-m^2$. In our setting it corresponds to $\Delta=4D>0$ and $a=-1<0$. In \cite{GrMo} they are considered two families of polynomials with $\Delta<0$ and $a>0$. To interpret correctly this claims, it is important to keep in mind that the exceptional nature of a discriminant depends on our scale and in some sense an exceptional discriminant, zero or character is like a sequence. The existence of a real zero in $[1-c/\log q,1]$ only ruins the generic de la Vall\'ee Poussin zero free region if $c$ can be taken arbitrarily small when $q$ grows. In \cite{GrMo} and \cite{FrIw} a bona fide asymptotic formula is only achieved if $\Delta$ is allowed to grow. 

\medskip

Our goal in this paper is to adapt the techniques of \cite{FrIw} to get a result valid for every $f$ as above when $\Delta$ is exceptional. 
By the reasons explained before we prefer to present the result as a main term plus an error term instead of as an asymptotic formula resembling the original statement of the conjectures.

\medskip

For each $N\ge 1$ we denote
\[
 \pi_f(N)=\#\{0\le f(n)\le N\,:\,f(n)\text{ prime}\}
\]
with $f$ as before 
and for each integer $d$, $\rho(d)=\#\{n\,:\,f(n)\equiv 0\pmod d\}$. Let   $\chi_\Delta(n)$  be the Kronecker  character modulo $\Delta$, $\chi_\Delta(n)=\left(\frac{\Delta}{n}\right)$. We will also consider the $L$ function associated to the character  
\begin{equation}\label{eq:ele}
L(s,\chi)= \sum_{n}\frac{\chi_{\Delta}(n)}{n^s}=\prod_{p\text{ prime}}\left(1-\frac{\chi_{\Delta}(p)}{p^s}\right)^{-1}.
\end{equation}
We will denote $\mathcal A=\big\{f(n)\in [0,N]\,:\, n\in\mathbb{Z} \big\}$ and $\mathcal A_d=\{k\in\mathcal A\,:\, d\mid k\}$   for $N,d\in\mathbb{Z}^+$ and  $A$ and $A_d$ stand for their cardinality. Also, we denote $V(x)=\prod_{p<x}\left(1-p^{-1}{\rho(p)}\right)$. The  exceptionality of the character $\chi_\Delta$ will be measured by 
\[
\beta=-\log(L(1,\chi_\Delta)\log|\Delta|).
\]
In the development of the proof  it is convenient  to introduce also 
\[
 L=-\log\left(L(1,\chi_{\Delta})\log A\right)\qquad \text{and}
 \qquad
 B=\frac{3\log{|\Delta|}}{\log A}.
\] 
To state our main result we introduce the function 
\[
g(\Delta)=\begin{cases}\Delta e^{-\beta/2} &\text{ if } \Delta>0\\ 
{|\Delta|}\big({4|a|- e^{-\beta/2}}\big)^{-1} &\text{ if } \Delta<0.
\end{cases}
\]
\begin{thm}\label{teo:main}
Let $1\le |a|\le e^{\beta/5}$ and   $g(\Delta)\le N\le  |a||\Delta|^{\beta/2}$.  Then
\[
\pi_f(N)=AV(A)\big(1+O(e^{-\sqrt{\beta}/6})\big)
\]
with an absolute $O$-constant. 
\end{thm}

\begin{rem}
 An asymptotic formula is obtained only under $\beta\to+\infty$ or equivalently, under the usual definition of exceptionality   $L(1,\chi_\Delta)\log |\Delta|\to 0$.
\end{rem}

\section{Guidelines. } 

Along the proof we follow the sieve techniques of \cite{FrIw}.  As usual for  $z\ge 2$ we write $P(z)=\prod_{p< z}p$ and
$S(\mathcal A,z)=\#\{a\in\mathcal A\,:\, (a,P(z))=1\}$. We start with the trivial identity 
\[
\pi_f(N)=S(\mathcal A,\sqrt N)+O(1)
\]
and use the well known Buchstab identity
\begin{equation}\label{eq:buch}
S(\mathcal A,z)-S(\mathcal A,\sqrt N)=\sum_{z\le p\le \sqrt N}S(\mathcal A_p,p).
\end{equation}
The main term of the theorem will come from  $S(\mathcal A,z)$, while the sum on the right will be part of the error term. In any event, in order to estimate both terms we need to have a concrete knowledge of both $A$ and $A_p$.
Then, the whole idea of the proof of Theorem \ref{teo:main} is to bound  the right hand side of (\ref{eq:buch}) using the exceptionality of the character $\chi_\Delta$.  This comes by noting that for~$d$ squarefree $\rho(d)\le \lambda(d)$ where $\lambda(n)$ is given by the convolution  $\lambda=1*\chi_\Delta$, and in particular $S(\mathcal A_p,p)=0$ if $\chi_\Delta(p)=-1$. If the character is exceptional, this will happen often, giving many zero terms in the right hand side of  (\ref{eq:buch}).

In order to estimate the right hand side  clearly we will need  to have some control  over $A_p$.   Observe that
\[
A_d=\frac{\rho(d)}{d}A+r_d, \quad \text{ with } |r_d|\le \rho(d). 
\]
Thanks to $\rho(d)\le \lambda(d)$   finding good bounds for the sum  in the right hand side of (\ref{eq:buch}) will come from finding  good estimates for the sum defined as 
\begin{equation}\label{eq:delta}\delta(x)=\sum_{x\le p\le A}\frac{\lambda(p)}{p}.
\end{equation}

\section{Proof of Theorem \ref{teo:main}. } 

Along the proof we will assume $\beta$ large enough, since otherwise Theorem \ref{teo:main} is the classical upper bound from linear sieve theory. 
The size of $A$ grows with $\Delta$ (see Lemma~\ref{lem:sizeA} below) then we can assume that $A$ is bigger than a large constant. 
We will frequently use the inequality $\beta\le \varepsilon\log|\Delta|$, which follows from Siegel's theorem.

\medskip

Let us start by finding the asymptotics of $S(\mathcal A,z)$. For that we use the fundamental lemma of sieve theory (see e.g. \cite[Cor. 6.10]{opera}), with level of distribution $A^{2/3}$,  $z^s=A^{2/3}$ and any $1<s<\frac{2\log A}{9\log\log A}$   to get
\begin{equation}\label{eq:az}
S(A,z)=AV(z)(1+O(e^{-s}))+O(A^{2/3}\log A).
\end{equation}
Observe that
\[
V(z)=V(A)e^{O(\delta(z))}=V(A)\big(1+O(\delta(z))\big),
\]
and hence we can replace $V(z)$ by $V(A)$  with an  error term  bounded by $\delta(z)$, which will be absorbed in the error term in (\ref{eq:az}).

\medskip

The rest of the paper will be dedicated to  bound  the right hand side of (\ref{eq:buch}), which we will split into three different sums, depending on the range of summation for the primes.
\begin{eqnarray}
\sum_{z\le p\le \sqrt N}S(\mathcal A_p,p)&=&\sum_{z\le p\le A/z^2}
+\sum_{A/z^2< p\le A}
+\sum_{A< p\le \sqrt N}\nonumber \\
&=&S_1+S_2+S_3.\label{S123}
\end{eqnarray}

We start with the sum $S_1$. The trivial bound  $S(\mathcal A_p,p)\le A_p$ is not good enough for the small  primes in the sum $S_1$, and we need a better bound gotten, as in \cite{FrIw}, using an upper bound sieve of dimension $2$ and level of distribution $A/pz$ (see e.g. \cite[Cor. 6.10]{opera}). This gives us 
\[
S(\mathcal A_p, p)\ll \frac{\rho(p)}{p}AV(z)+\sum_{d<A/pz}\rho(d),
\]
where the sum runs over squarefree integers $d$.  The last term is trivially bounded by
\[
\sum_{d<A/pz}\rho(d)\le \sum_{d<A/pz}(1*\chi_\Delta)(d)\ll \sum_{d<A/pz}\tau(d)\ll \frac{A}{pz}\log A.
\]
Noting that  $z\ge (\log A)^3$, which follows by our assumption in $s$, and  that  $V(z)\ge V(A)\gg (\log A)^{-2}$, since $\rho(p)\le 2$ for any prime $p$, we end up with
$
S(\mathcal A_p, p)\ll p^{-1}{\lambda(p)}AV(z)
$
and hence
\begin{equation}\label{eq:aps}
S_1\ll AV(z)\delta(z).
\end{equation}

\

To bound  $\delta(z)$ we use Lemma 3.4 of \cite{FrIw}, which we include for reader's convenience.
\begin{lem}\label{lem:friwi}  Let $2\le u\le y\le x$. Then, 
\[
\sum_{\substack{y\le n\le x\\(n,P(u))=1}}\frac{\lambda(n)}{n}\ll W(u)L(1,\chi_\Delta)\log\left(\frac xy\right)
+|\Delta|^{1/8+\varepsilon}y^{-1/3}u^{1/3}\log u,
\]
where $W(u)=\prod_{p<u}\left(1-p^{-1}\right)\left(1-p^{-1}{\chi_\Delta(p)}\right)$.
\end{lem}

\begin{rem}
 Observe that, since  
$W(u)\le C\prod_{p<u}\left(1-p^{-1}{\lambda(p)}\right)\le CV(u)$, for some absolute constant $C$, we can write either $W(u)$ or $V(u)$ indistinctly. 
\end{rem}

Further we will  use the formula,  also proved in \cite[p.1106]{FrIw},
\begin{equation}\label{eq:delfi}
\delta(z)^k\ll kk!W(z)L(1,\chi_{\Delta})\log A
+k!|\Delta|^{1/8+\varepsilon}z^{(1-k)/3}
\end{equation}
valid for any integer $k\ge 1$.  It is worth to note that in order to establish the previous formula it is needed a bound of the type $\log z\ll \Delta^{\epsilon}$, which in our case follows assuming  $L>0$. Indeed
\[
\log z<\log A<\frac1{L(1,\chi_\Delta)}\ll\Delta^{\epsilon}.
\]
Dropping the contribution of $W(z)$ in (\ref{eq:delfi}) we get the more convenient form
\[
\delta(z)\ll k\left(kL(1,\chi_{\Delta})\log A
+|\Delta|^{1/8+\varepsilon}z^{(1-k)/3}\right)^{1/k}.
\]

\

Our goal is to prove 
\begin{equation}\label{eq:error}
\delta(z)\ll e^{-s}
\end{equation}
with a proper selection of $s$. 

Taking any positive integer $k\ge(\frac3{16}+3\varepsilon)Bs+1$ and noting that $z=|\Delta|^{2/Bs}$ we obtain 
\[
|\Delta|^{1/8+\varepsilon}z^{(1-k)/3}\ll |\Delta|^{-\epsilon}.
\]
On the other hand
\[
kL(1,\chi_{\Delta})\log A\ge |\Delta|^{-\epsilon}
\]
follows by Siegel's theorem, and then
\[
\delta(z)\ll k(L(1,\chi_{\Delta})\log A)^{1/k}.
\]
Observe that, assuming again $L>0$,  we have that the previous bound is increasing in $k$, and so we can relax the condition of $k$ being an integer, In particular, we can  take $k=Bs$, which is possible assuming $Bs$ greater than a constant greater than ${16}/{13}$. Then, to prove (\ref{eq:error}) we need to select some 
\[
s\le -\log k+\frac{L}{k},
\]
which gives, replacing the value of  $k$,
\begin{equation}
Bs^2+Bs\log(Bs)\le L.
\end{equation}
The error term in Theorem \ref{teo:main} is in terms of $\beta$ instead of $L$. The comparison between both quantities comes from a proper control in $A$. We have the following  lemma. 

\begin{lem}\label{lem:sizeA} Let $A\ge 1$ and assume the hypothesis in Theorem \ref{teo:main}. Given $\varepsilon>0$, we have
\[
 |\Delta|^{1/2-\varepsilon}\le A< \frac{4\sqrt N}{\sqrt{|a|}},
\]
for $\Delta$ large enough (depending on $\varepsilon$).
\end{lem}

\begin{rem}
 For the application of this in the proof of the main result we are going to choose $1/2-\varepsilon=7/16$. This is connected to the constant  $16/13$ above. 
\end{rem}

\begin{proof}
 The inequalities $0\le f(x)\le N$ define one or two intervals for $x$, depending on the real zeros of  $f$ and the sign of $a$ and $\Delta$, and it is straighforward to measure the length of those intervals to be
\[
X= \begin{cases} 
\frac{\sqrt {\Delta+4aN}}{a} & \text{ if } \Delta<0,\,a>0,\,N> \frac{|\Delta|}{4|a|},\\
\frac{4N}{\sqrt {\Delta+4aN}+\sqrt{\Delta}} & \text{ if } \Delta>0,\,a>0\text{ or }\text{ if } \Delta>0,\,a<0,\,N\le \frac{|\Delta|}{4|a|}, \\
\frac{\sqrt {\Delta}}{|a|} & \text{ if } \Delta>0,a<0,N> \frac{|\Delta|}{4|a|}.
\end{cases}
\]
The cases not listed above give empty intervals.  
From here the upper bound $X\le 2\sqrt{N/{|a|}}$ is trivial. Then, noting 
\begin{equation}\label{eq:ax}
X-2<A<X+2,
\end{equation}
we deduce 
$A\le 4\sqrt{N/|a|}$ for  $|a|<N$, whenever $X\ge 2$, which follows  from our assumptions $a<e^{\beta/5}$, $N>g(\Delta)$ by Siegel's theorem because $\beta$ and $\Delta$ can be assumed sufficiently large.
 
We now prove $X>|\Delta|^{1/2-\varepsilon}$. In  the last case in the definition of $X$, the result follows again from $|a|<e^{\beta/5}$. If $\Delta<0$, then $X>a^{-1}{e^{-\beta/4}}\sqrt N$, since $N\ge g(\Delta)$. Further if $\Delta>0$ and $N\le \frac{\Delta}{4|a|}$, we have $X\gg {N}/{\sqrt\Delta}$, finally if  $\Delta>0$, $a>0$ and $N> \frac{\Delta}{4|a|}$ we have the stronger bound $X\gg\sqrt{ N/|a|}$. In any case $X>\Delta^{1/2-\varepsilon}$ 
is a consequence of $N/\Delta\gg |\Delta|^{-\epsilon}$ and $a\ll |\Delta|^\epsilon$.
\end{proof}
 
Now, $A\le\frac{4\sqrt N}{\sqrt {|a|}}$ and the upper bound for  $N$ give
\[
\beta=L+\log\left(\frac{\log A}{\log |\Delta|}\right)<2L.
\]

We select $s=\frac 12 \sqrt{{\beta}/{B}}$.  
By Lemma~\ref{lem:sizeA} $B$ is bounded, namely with the choice of $\varepsilon$ as in the Remark we have $B<7$. Then $s$ is arbitrarily large, in particular $s>1$. 
It is important to check that this selection of $s$ is compatible with the rest of our previous assumptions: 
\[
\frac{16}{13B}<s<\frac{2\log A}{9\log\log A}.
\] 
The first inequality is consequence of $A<4\sqrt{N/|a|}$ and the upper bound in $N$. The  second is equivalent to 
\[
1<\frac{4\sqrt B\log A}{9\sqrt\beta\log\log A}, 
\]
which follows for $A$ large enough, by the definition of $B$ and Siegel's Theorem.

Let us prove with this selection of $s$ that 
\begin{equation}\label{eq:ese}
Bs^2+Bs\log(Bs)\le \frac\beta2<L.
\end{equation}
As $B<7$ and $s$ is arbitrarily large, we can suppose  $Bs^2\ge Bs\log(Bs)$ and then (\ref{eq:ese}) follows directly from our choice of $s$. 
This proves (\ref{eq:ese}), and (\ref{eq:error}) with $s=\frac 12 \sqrt{{\beta}/{B}}$ (and assures $L>0$ as assumed), which gives by~\eqref{S123}
\begin{equation}\label{S1b}
 S_1\ll AV(z)e^{-\sqrt{\beta/4B}}\ll AV(A)e^{-\sqrt{\beta}/6}
\end{equation}
since  $B<7$.

\medskip

It remains to bound $S(\mathcal A_p, p)$ for medium and large $p$. If $A/z^2<p<A$ then
\[
S(\mathcal A_p,p)\le A_p\le \frac{\rho(p)}{p}A+\lambda(p)\ll \frac{\lambda(p)}{p}A,
\]
so
\begin{equation}\label{eq:apsbig}
S_2\ll A \delta(A/z^2).
\end{equation}
We apply  Lemma \ref{lem:friwi}  $x=A$ and $y=A/z^2$ and $u=Az^{-2}|\Delta|^{-3/8-4\varepsilon}$.  Observe that $u=A^\gamma$
with 
$\gamma=1-\frac4{3s}-\frac B8-\frac43 B\varepsilon$ and $\gamma>0$ 
for $\varepsilon$ sufficiently small and $s$ sufficiently large,  since $B<7$. With this selection 
the first term in the sum in  Lemma~\ref{lem:friwi} dominates the second and we deduce
\[
\delta(A/z^2)\ll V(u)L(1,\chi_\Delta)\log(z^2)\ll V(A) L(1,\chi_\Delta)\log(z^2).
\]
For the last inequality we have used that $u$ is a positive power of $A$. Now,  
\[
\log(z^2)L(1,\chi_\Delta)=\frac4{3s}e^{-L}<e^{-\beta/2}\le e^{-\sqrt{\beta}/6},
\]
and putting everything together we get the desired result
\begin{equation}\label{S2b}
 S_2\ll AV(A)e^{-\sqrt{\beta}/6}.
\end{equation}

Finally, it remains to bound $S_3$ corresponding to the primes $A<p<\sqrt N$. We have 
\[
S(\mathcal A_p,p)\le A_p\le \frac{\rho(p)}{p}A+\lambda(p)\ll \frac{\lambda(p)}{p}\sqrt N,
\]
and again using Lemma \ref{lem:friwi} with the same parameter $u$ as before, $x=\sqrt N$ and $y=A$, we get
\[
\sum_{A\le p\le \sqrt N}S(\mathcal A_p,p)\ll \sqrt N V(A)L(1,\chi_{\Delta})\log\sqrt N
\]
We separate the different cases giving $f(\mathbb{Z})\cap \mathbb{Z}^+\ne\emptyset$ (see the definition of $X$ in the proof of Lemma~\ref{lem:sizeA}). 

If $\Delta<0$,  $a>0$ and  by hypothesis $N \ge g(\Delta)$, which  gives,  $\sqrt N\ll ae^{\beta/4}A$ by (\ref{eq:ax}) and, hence, since  $f(x)=x\log x$ is increasing for $x>2$, we get 
\[
\sqrt N \log\sqrt N\ll ae^{\beta/4}A\log(ae^{\beta/ 4}A)\ll ae^{\beta/ 4}A\log A
\]
by our assumption in $a$ and Lemma \ref{lem:sizeA}. Also $a<e^{\beta/5}$ implies
\[
a\ll e^{L-\beta/4-{\sqrt \beta}/{6}},
\]
since $\beta<2L$, which gives
\begin{equation}\label{S3b}
 S_3\ll AV(A)e^{-{\sqrt\beta}/{6}},
\end{equation}
in this case as desired.

If $\Delta>0$, $a<0$, we just need to consider the case $\Delta e^{-\beta/2}\le N\le \frac{\Delta}{4|a|}$, since  $f(n)\le  \frac{\Delta}{4|a|}$. Then  $A\gg{N}/{\sqrt\Delta}>e^{-\beta/4}\sqrt N$ and the proof of \eqref{S3b} follows in the same way as before. 

Finally, if $\Delta>0$, $a>0$ we have $A\gg \sqrt{N/a}\gg e^{-\beta/10}\sqrt N$ and the same proof applies getting again \eqref{S3b}. 

\medskip

Substituting \eqref{S1b}, \eqref{S2b} and \eqref{S3b} in \eqref{S123} and recalling \eqref{eq:buch} and \eqref{eq:az}, the proof of Theorem~\ref{teo:main} is complete.


\begin{thebibliography}{1}

\bibitem{opera}
J.~Friedlander and H.~Iwaniec.
\newblock {\em Opera de cribro}, volume~57 of {\em American Mathematical
  Society Colloquium Publications}.
\newblock American Mathematical Society, Providence, RI, 2010.

\bibitem{FrIw}
J.~B. Friedlander and H.~Iwaniec.
\newblock Exceptional discriminants are the sum of a square and a prime.
\newblock {\em Q. J. Math.}, 64(4):1099--1107, 2013.

\bibitem{GrMo}
A.~Granville and R.~A. Mollin.
\newblock Rabinowitsch revisited.
\newblock {\em Acta Arith.}, 96(2):139--153, 2000.

\bibitem{HaLi}
G.~H. Hardy and J.~E. Littlewood.
\newblock Some problems of `{P}artitio numerorum'; {III}: {O}n the expression
  of a number as a sum of primes.
\newblock {\em Acta Math.}, 44(1):1--70, 1923.

\end{thebibliography}

\end{document}